\documentclass{notemat}
\usepackage[utf8]{inputenc}
\usepackage[T1]{fontenc}
\usepackage{amsmath}
\usepackage{amssymb}
\usepackage{xcolor}
\usepackage{nmmacro}

\newcommand{\coloneqq}[0]{\mathrel{\mathop:}=}
\newcommand{\notproves}[2]{{#1} \not\vdash {#2}}
\newcommand{\Sigmafano}[0]{\Sigma^{\prime}}
\newcommand{\Sigmasharpened}[0]{\Sigma_{*}}
\newcommand{\Aaxiom}[1]{\mathrm{A}{#1}}
\newcommand{\Baxiom}[1]{\mathrm{B}{#1}}
\newcommand{\sigmafunc}[2]{{\sigmasymbol}_{#1}({#2})}
\newcommand{\sigmasymbol}[0]{\sigma}
\newcommand{\tausymbol}[0]{\tau}
\newcommand{\taufunc}[3]{{\tausymbol}_{{#1}{#2}}({#3})}
\newcommand{\collinearsymbol}[0]{L}
\newcommand{\collinear}[3]{{\collinearsymbol}({#1}{#2}{#3})}
\newcommand{\paradox}[0]{\textsf{Paradox}}
\newcommand{\macefour}[0]{\textsf{Mace4}}
\newcommand{\tipi}[0]{\textsf{Tipi}}
\newcommand{\subtract}[2]{{#1} \setminus \{ {#2} \}}
\newcommand{\add}[2]{{#1} \cup \{ {#2} \}}
\newcommand{\proves}[2]{{#1} \vdash {#2}}

\let\phi=\varphi
\begin{document}
\begin{article}
\begin{opening}
\title{Complete independence of an axiom system for central translations}
\author{Jesse Alama}
\institute{Center for Artificial Intelligence, New University of Lisbon, Portugal\email{j.alama@fct.unl.pt}}
\runningauthor{Alama}
\runningtitle{Complete independence for central translations}
\begin{abstract}A recently proposed axiom system for André's central translation structures is improved upon.  First, one of its axioms turns out to be dependent (derivable from the other axioms).  Without this axiom, the axiom system is indeed independent.  Second, whereas most of the original independence models were infinite, finite independence models are available.  Moreover, the independence proof for one of the axioms employed proof-theoretic techniques rather than independence models; for this axiom, too, a finite independence model exists.  For every axiom, then, there is a finite independence model.  Finally, the axiom system (without its single dependent axiom) is not only independent, but completely independent.\end{abstract}
\end{opening}
\section{Introduction}

\par{Pambuccian has offered two axiom systems $\Sigma$ and $\Sigmafano$~\cite{pambuccian2001two} for André's central translation structures~\cite{andre1961parallelstrukturen}.  ($\Sigmafano$ is $\Sigma$ together with the Fano principle that diagonals of parallelograms are not parallel.) In this note we further develop Pambuccian's work by showing that:}

\begin{itemize}
\item {$\Sigma^{\prime}$ has a dependent axiom.  Without this axiom, the axiom system is indeed independent.}
\item {Without the dependent axiom, finite independence models exist for all axioms, whereas most of the independence proofs offered in~\cite{pambuccian2001two} used infinite models.  In particular, one of the independence proofs was accomplished by proof-theoretic methods rather than by a independence model; but for this axiom, too, there is a finite independence model;}
\item {Without its dependent axiom, $\Sigma^{\prime}$ is not just independent, but completely independent.}
\end{itemize}

\par{For the sake of completeness, we repeat here the definitions of the axiom systems $\Sigma$ and $\Sigmafano$.  Both are based on classical one-sorted first-order logic with identity.  Variables are intended to range over points.  A single language is used for both axiom systems; the language has three constants, $a_{0}$, $a_{1}$, and $a_{2}$, a ternary relation $L$ for collinearity, and a single ternary function symbol $\tau$ for central translations: $\taufunc{a}{b}{c}$ is the image of $c$ under the translation that shifts $a$ to $b$.  With the the binary operation symbol $\sigma$ understood as}
\[
\sigmafunc{a}{b} = \taufunc{b}{a}{a},
\]
\par{(essentially a point reflection), the axioms of $\Sigma$ are as follows:}

\begin{description}
\item[A3] $a \ne b \wedge \collinear{a}{b}{c} \wedge \collinear{a}{b}{d} \rightarrow \collinear{a}{c}{d}$
\item[B1] $\collinear{a}{b}{c} \rightarrow \collinear{b}{a}{c}$
\item[B2] $\taufunc{a}{b}{c} = \taufunc{a}{c}{b}$
\item[B3] $\collinear{a}{b}{\sigmafunc{a}{b}}$
\item[B4] $\collinear{a}{b}{c} \rightarrow \collinear{x}{\taufunc{a}{b}{x}}{\taufunc{a}{c}{x}}$
\item[B5] $\taufunc{a}{b}{x} = x \rightarrow a = b$
\item[B6] $\taufunc{a}{b}{x} = \taufunc{c}{\taufunc{a}{b}{x}}{x}$
\item[B7] $\neg \collinear{a_{0}}{a_{1}}{a_{2}}$
\end{description}
\par{(The appearance of $\Aaxiom{3}$ without $\Aaxiom{1}$ and $\Aaxiom{2}$ is not an error.  In the the official definition of $\Sigma$ from~\cite{pambuccian2001two}, rather than duplicating an axiom from André's axiom system, the names of whose axioms all have the prefix ``A'', under a new name, it is simply reused and the ``B'' axioms are offered.)  The axiom system $\Sigma^{\prime}$ is $\Sigma$ together with}
\begin{description}
\item[B8] $\sigmafunc{a}{b} = b \rightarrow a = b$
\end{description}
\par{$\Baxiom{8}$ captures the Fano principle that diagonals of parallelograms are not parallel.}
\par{Universal quantifiers will usually be suppressed.}
\section{A dependent axiom}
\par{It was claimed in~\cite{pambuccian2001two} that $\Sigmafano$ is independent.  This requires qualification: $\Baxiom{5}$ is indeed an independent axiom of $\Sigma$ (Proposition~\ref{b5-independent-of-sigma}), but $\Baxiom{5}$ can be proved from $\subtract{\Sigma}{\Baxiom{5}}$ with the help of the Fano principle $\Baxiom{8}$ (Proposition~\ref{b5-dependent}).}
\begin{Proposition}\label{b5-independent-of-sigma}$\notproves{\subtract{\Sigma}{\Baxiom{5}}}{\Baxiom{5}}$.\end{Proposition}
\begin{proof}
\par{Consider the domain $\{1,2\}$ and interpret $\collinearsymbol$ and $\tausymbol$ according to Table~\ref{tab:b5-counterexample}. A counterexample to}
\[
\taufunc{a}{b}{x} = x \rightarrow a = b
\]
\par{is provided by $(a,b,x) \coloneqq (1,2,2)$.}
\begin{table}
\begin{tabular}{lll}
\hline
Triple & $\tausymbol$ & $\collinearsymbol$\\
\hline
$(1,1,1)$ & $2$ & $-$\\
$(1,1,2)$ & $2$ & $+$\\
$(1,2,1)$ & $2$ & $-$\\
$(1,2,2)$ & $2$ & $+$\\
$(2,1,1)$ & $2$ & $-$\\
$(2,1,2)$ & $2$ & $+$\\
$(2,2,1)$ & $2$ & $-$\\
$(2,2,2)$ & $2$ & $-$\\
\end{tabular}\caption{\label{tab:b5-counterexample}A model of $\add{\subtract{\Sigma}{\Baxiom{5}}}{\neg\Baxiom{5}}$}
\end{table}
For each pair $(a,b)$, ${\tausymbol}_{a,b}$ fails to be a transitive action because $1$ is never a value of $\tausymbol$.
\end{proof}
%% \begin{Lemma}\label{b5-dependent-lemma-1}$\proves{\Sigmafano}{\taufunc{a}{a}{a} = b \rightarrow a = b}$
%% \end{Lemma}
%% \begin{proof}
%% \par{The conclusion follows from}
%% \[
%% \sigmafunc{a}{a} = a,
%% \]
%% \par{which is equation~(2) of~\cite{pambuccian2001two}.  This result does not depend on $\Baxiom{5}$; in fact, $\{ \Baxiom{2}, \Baxiom{6}, \Baxiom{8} \}$ suffices.}
%% \end{proof}
%% \par{Lemma~\ref{b5-dependent-lemma-1} fails in the absence of $\Baxiom{8}$.}
\begin{Lemma}\label{b5-dependent-lemma-2}$\proves{\Sigmafano}{\taufunc{a}{b}{c} = \taufunc{d}{c}{\taufunc{a}{b}{d}}}$
\end{Lemma}
\begin{proof}{The desired conclusion follows from $\Baxiom{2}$ and $\Baxiom{6}$ ($\Baxiom{5}$ is not needed).}\end{proof}
\begin{Lemma}\label{b5-dependent-lemma-3}$\proves{\Sigmafano}{\taufunc{a}{b}{a} = b}$\end{Lemma}
\begin{proof}This is equation (3) of~\cite{pambuccian2001two}.  It is derived without the help of $\Baxiom{5}$.  (Indeed, $\{ \Baxiom{2}, \Baxiom{6}, \Baxiom{8} \}$ suffices).\end{proof}
%% \begin{Lemma}\label{b5-dependent-lemma-4}$\proves{\Sigmafano}{\taufunc{a}{b}{\taufunc{b}{a}{c}} = c}$\end{Lemma}
%% \begin{proof}
%% \end{proof}
%% \par{The proof of Lemma~\ref{b5-dependent-lemma-4} uses the Fano principle.}
\begin{Proposition}\label{b5-dependent}$\proves{\subtract{\Sigmafano}{\Baxiom{5}}}{\Baxiom{5}}$.\end{Proposition}
\begin{proof}
\par{Suppose $\taufunc{a}{b}{x} = x$.  From Lemmas~\ref{b5-dependent-lemma-2} and~\ref{b5-dependent-lemma-3}, we have}
\begin{equation}\label{eq:1}
\taufunc{u}{v}{w} = \taufunc{u}{w}{v},
\end{equation}
\par{as well as}
\begin{equation}\label{eq:2}
\taufunc{\taufunc{t}{u}{v}}{w}{u}  = \taufunc{v}{t}{w}
\end{equation}
\par{for all $t$, $u$, $v$, and $w$.  Thus, by~(\ref{eq:1}), $\taufunc{a}{x}{b} = x$, whence (by Lemma~\ref{b5-dependent-lemma-2} and~(\ref{eq:2}),}
\begin{equation*}
\taufunc{u}{a}{v} = \taufunc{u}{b}{v}
\end{equation*}
\par{for all $u$ and $v$.  Lemma~\ref{b5-dependent-lemma-3} then gives us the desired conclusion $a = b$.}
\end{proof}
\par{In light of Proposition~\ref{b5-dependent}, we define:}
\begin{Definition}$\Sigmasharpened \coloneqq \Sigmafano \setminus \{ \Baxiom{5} \}$.\end{Definition}
\par{In the following, the theory in focus is $\Sigmasharpened$ rather than $\Sigmafano$.}
\section{Small finite independence models}\label{sec:finite-independence-models}
\par{The next several propositions show that every axiom of $\Sigmasharpened$ has a finite independence model.  (Incidentally, the cardinalities of the independence models are minimal: when it is claimed that there is a independence model for $\phi$ of cardinality $n$, it is also claimed that $\phi$ is true in every model of $\subtract{\Sigmasharpened}{\phi}$ of cardinality less than $n$.)}
\par{In the independence models that follow we give only the interpretation of the predicate $\collinearsymbol$ and the function $\tausymbol$.  Strictly speaking, this is not enough, because we need to interpret the constants $a_{0}$, $a_{1}$, and $a_{2}$ so that $\Baxiom{7}$ holds.  But in the independence models there is always at least one triangle; from any one, an interpretation of the constants $a_{0}$, $a_{1}$, and $a_{2}$ can be chosen so that $\Baxiom{7}$ is satisfied.}
\begin{Proposition}\label{a3-finite-model}There exists an independence model for $\Aaxiom{3}$ of cardinality~$3$.\end{Proposition}
\begin{proof}
\par{Without $\Aaxiom{3}$ one cannot prove}
\[
\collinear{a}{a}{b} \wedge \collinear{a}{b}{a} \wedge \collinear{a}{b}{b},
\]
\par{the failure of which opens the door to geometrically counterintuitive models.  Consider the domain $\{1,2,3\}$, interpret $\collinearsymbol$ and $\tausymbol$ as in Table~\ref{tab:a3-counterexample}.  Note that the model is ``collinear'' in the sense that $\collinear{\alpha(1)}{\alpha(2)}{\alpha(3)}$ for any permutation $\alpha$ of $\{ 1, 2 3 \}$; nonetheless, for many pairs $(u,v)$ of distinct points, the various collinearity statements one can make about $u$ and $v$ are false.  A counterexample to}
\[
a \ne b \wedge \collinear{a}{b}{c} \wedge \collinear{a}{b}{d} \rightarrow \collinear{a}{c}{d}
\]
\par{is $(a,b,c,d) \coloneqq (3,2,1,1)$.}
\begin{table}[h]
\begin{tabular}{lll|lll|lll}
\hline
Triple & $\tausymbol$ & $\collinearsymbol$ & Triple & $\tausymbol$ & $\collinearsymbol$ & Triple & $\tausymbol$ & $\collinearsymbol$\\
\hline
$(1,1,1)$ & $1$ & $+$ & $(2,1,1)$ & $3$ & $+$ & $(3,1,1)$ & $2$ & $-$\\
$(1,1,2)$ & $2$ & $-$ & $(2,1,2)$ & $1$ & $-$ & $(3,1,2)$ & $3$ & $+$\\
$(1,1,3)$ & $3$ & $+$ & $(2,1,3)$ & $2$ & $+$ & $(3,1,3)$ & $1$ & $+$\\
$(1,2,1)$ & $2$ & $+$ & $(2,2,1)$ & $1$ & $+$ & $(3,2,1)$ & $3$ & $+$\\
$(1,2,2)$ & $3$ & $-$ & $(2,2,2)$ & $2$ & $+$ & $(3,2,2)$ & $1$ & $+$\\
$(1,2,3)$ & $1$ & $+$ & $(2,2,3)$ & $3$ & $-$ & $(3,2,3)$ & $2$ & $-$\\
$(1,3,1)$ & $3$ & $-$ & $(2,3,1)$ & $2$ & $+$ & $(3,3,1)$ & $1$ & $-$\\
$(1,3,2)$ & $1$ & $+$ & $(2,3,2)$ & $3$ & $+$ & $(3,3,2)$ & $2$ & $+$\\
$(1,3,3)$ & $2$ & $+$ & $(2,3,3)$ & $1$ & $-$ & $(3,3,3)$ & $3$ & $+$\\
\end{tabular}\caption{\label{tab:a3-counterexample}A model of $\add{\subtract{\Sigmasharpened}{\Aaxiom{3}}}{\neg\Aaxiom{3}}$}
\end{table}
\end{proof}
\begin{Proposition}\label{b1-finite-model}There exists an independence model for $\Baxiom{1}$ of cardinality~$3$.\end{Proposition}
\begin{proof}
\par{The difficulty here is that}
\[
\collinear{a}{b}{a}
\]
\par{fails without $\Baxiom{1}$.  Consider the domain $\{ 1,2,3 \}$ and the interpretations of $\collinearsymbol$ and $\tausymbol$ as in Table~\ref{tab:b1-counterexample}.  From the standpoint of $\collinearsymbol$, the model is nearly trivialized; the only triples $(a,b,c)$ where $\collinear{a}{b}{c}$ fails are where $a \ne b$ and $a = c$.  An example where}
\[
\collinear{a}{b}{c} \rightarrow \collinear{b}{a}{c}
\]
\par{fails is $(a,b,c) \coloneqq (1,2,3)$.}
\begin{table}
\begin{tabular}{lll|lll|lll}
\hline
Triple & $\tausymbol$ & $\collinearsymbol$ & Triple & $\tausymbol$ & $\collinearsymbol$ & Triple & $\tausymbol$ & $\collinearsymbol$\\
\hline
$(1,1,1)$ & $1$ & $+$ & $(2,1,1)$ & $3$ & $+$ & $(3,1,1)$ & $2$ & $+$\\
$(1,1,2)$ & $2$ & $+$ & $(2,1,2)$ & $1$ & $-$ & $(3,1,2)$ & $3$ & $+$\\
$(1,1,3)$ & $3$ & $+$ & $(2,1,3)$ & $2$ & $+$ & $(3,1,3)$ & $1$ & $-$\\
$(1,2,1)$ & $2$ & $-$ & $(2,2,1)$ & $1$ & $+$ & $(3,2,1)$ & $3$ & $+$\\
$(1,2,2)$ & $3$ & $+$ & $(2,2,2)$ & $2$ & $+$ & $(3,2,2)$ & $1$ & $+$\\
$(1,2,3)$ & $1$ & $+$ & $(2,2,3)$ & $3$ & $+$ & $(3,2,3)$ & $2$ & $-$\\
$(1,3,1)$ & $3$ & $-$ & $(2,3,1)$ & $2$ & $+$ & $(3,3,1)$ & $1$ & $+$\\
$(1,3,2)$ & $1$ & $+$ & $(2,3,2)$ & $3$ & $-$ & $(3,3,2)$ & $2$ & $+$\\
$(1,3,3)$ & $2$ & $+$ & $(2,3,3)$ & $1$ & $+$ & $(3,3,3)$ & $3$ & $+$\\
\end{tabular}\caption{\label{tab:b1-counterexample}A model of $\add{\subtract{\Sigmasharpened}{\Baxiom{1}}}{\neg\Baxiom{1}}$}
\end{table}
\end{proof}

\begin{Proposition}\label{b2-finite-model}There exists an independence model for $\Baxiom{2}$ of cardinality~$3$.\end{Proposition}
\begin{proof}
\par{Consider the domain $\{ 1,2,3 \}$ and the interpretations of $\collinearsymbol$ and $\tausymbol$ are as in Table~\ref{tab:b2-counterexample}. A counterexample to}
\[
\taufunc{a}{b}{c} = \taufunc{a}{c}{b}
\]
\par{in this structure is $(a,b,c) \coloneqq (1,3,1)$.}
\begin{table}
\begin{tabular}{lll|lll|lll}
\hline
Triple & $\tausymbol$ & $\collinearsymbol$ & Triple & $\tausymbol$ & $\collinearsymbol$ & Triple & $\tausymbol$ & $\collinearsymbol$\\
\hline
$(1,1,1)$ & $1$ & $+$ & $(2,1,1)$ & $1$ & $+$ & $(3,1,1)$ & $1$ & $+$\\
$(1,1,2)$ & $1$ & $+$ & $(2,1,2)$ & $1$ & $+$ & $(3,1,2)$ & $1$ & $+$\\
$(1,1,3)$ & $1$ & $-$ & $(2,1,3)$ & $1$ & $-$ & $(3,1,3)$ & $1$ & $-$\\
$(1,2,1)$ & $2$ & $+$ & $(2,2,1)$ & $2$ & $+$ & $(3,2,1)$ & $2$ & $+$\\
$(1,2,2)$ & $2$ & $+$ & $(2,2,2)$ & $2$ & $+$ & $(3,2,2)$ & $2$ & $+$\\
$(1,2,3)$ & $2$ & $-$ & $(2,2,3)$ & $2$ & $-$ & $(3,2,3)$ & $2$ & $-$\\
$(1,3,1)$ & $2$ & $+$ & $(2,3,1)$ & $1$ & $+$ & $(3,3,1)$ & $2$ & $+$\\
$(1,3,2)$ & $2$ & $+$ & $(2,3,2)$ & $1$ & $+$ & $(3,3,2)$ & $2$ & $+$\\
$(1,3,3)$ & $2$ & $-$ & $(2,3,3)$ & $1$ & $-$ & $(3,3,3)$ & $2$ & $-$\\
\end{tabular}\caption{\label{tab:b2-counterexample}A model of $\add{\subtract{\Sigmasharpened}{\Baxiom{2}}}{\neg\Baxiom{2}}$}
\end{table}
\end{proof}
\begin{Proposition}\label{b3-finite-model}There exists an independence model for $\Baxiom{3}$ of cardinality~$1$.\end{Proposition}
\begin{proof}
\par{$\Baxiom{3}$ is the only axiom that outright asserts that some points are collinear.  Thus, if every triple $(a,b,c)$ of points constitutes a triangle (that is, $\neg\collinear{a}{b}{c}$ holds for all $a$, $b$, and $c$), then clearly $\Baxiom{3}$ would be falsified.  So take a $1$-element structure and interpret $\collinearsymbol$ so that $\collinear{a}{b}{c}$ for the unique triple $(a,b,c)$ of the structure is false.  The interpretations of $\tau$, $a_{0}$, $a_{1}$, and $a_{2}$ are forced.  One can check that all axioms, except $\Baxiom{3}$, are satisfied.}
\end{proof}
\par{Pambuccian was unable to find an independence model of $\Baxiom{4}$.  To show that $\Baxiom{4}$ is independent of the other axioms, methods of structural proof analysis~\cite{negri2008structural} were employed.  Specifically, an analysis of all possible formal derivations starting from $\Sigma$ was made, and by a syntactic-combinatorial argument it was found that $\Baxiom{4}$ is underivable from $\subtract{\Sigma}{\Baxiom{4}}$.  By the completeness theorem, then, there must exist an independence model.  Here is one:}
\begin{Proposition}\label{b4-finite-model}There exists an independence model for $\Baxiom{4}$ of cardinality~$27$.\end{Proposition}
\begin{proof}
\par{For lack of space, we omit an explicit description of the $27 \times 27 \times 27$ table.  The model was found by the finite model-finding program {\macefour}~\cite{prover9-mace4}.}
\end{proof}
\par{Consideration of $\Baxiom{5}$ is skipped because it is a dependent axiom of $\Sigma$ (and in any case is not officially an axiom of $\Sigmasharpened$).}
\begin{Proposition}\label{b6-finite-model}There exists a independence model for $\Baxiom{6}$ of cardinality~$3$.\end{Proposition}
\begin{proof}
\par{Consider the domain $\{ 1,2,3 \}$ and the interpretations of $\collinearsymbol$ and $\tausymbol$ as in Table~\ref{tab:b6-counterexample}.  A counterexample to}
\[
\taufunc{a}{b}{x} = \taufunc{c}{\taufunc{a}{b}{x}}{x}
\]
is given by $(a,b,c,x) \coloneqq (2,2,1,2)$.
\begin{table}
\begin{tabular}{lll|lll|lll}
\hline
Triple & $\tausymbol$ & $\collinearsymbol$ & Triple & $\tausymbol$ & $\collinearsymbol$ & Triple & $\tausymbol$ & $\collinearsymbol$\\
\hline
$(1,1,1)$ & $2$ & $+$ & $(2,1,1)$ & $1$ & $+$ & $(3,1,1)$ & $2$ & $+$\\
$(1,1,2)$ & $2$ & $+$ & $(2,1,2)$ & $1$ & $+$ & $(3,1,2)$ & $2$ & $+$\\
$(1,1,3)$ & $2$ & $-$ & $(2,1,3)$ & $2$ & $-$ & $(3,1,3)$ & $2$ & $-$\\
$(1,2,1)$ & $2$ & $+$ & $(2,2,1)$ & $1$ & $+$ & $(3,2,1)$ & $2$ & $+$\\
$(1,2,2)$ & $2$ & $+$ & $(2,2,2)$ & $1$ & $+$ & $(3,2,2)$ & $2$ & $+$\\
$(1,2,3)$ & $2$ & $-$ & $(2,2,3)$ & $2$ & $-$ & $(3,2,3)$ & $2$ & $-$\\
$(1,3,1)$ & $2$ & $+$ & $(2,3,1)$ & $2$ & $+$ & $(3,3,1)$ & $2$ & $+$\\
$(1,3,2)$ & $2$ & $+$ & $(2,3,2)$ & $2$ & $+$ & $(3,3,2)$ & $2$ & $+$\\
$(1,3,3)$ & $2$ & $-$ & $(2,3,3)$ & $1$ & $-$ & $(3,3,3)$ & $2$ & $-$\\
\end{tabular}\caption{\label{tab:b6-counterexample}A model of $\add{\subtract{\Sigmasharpened}{\Baxiom{6}}}{\neg\Baxiom{6}}$}
\end{table}
\end{proof}
\par{We are unable to improve upon the independence proof for $\Baxiom{7}$ given in~\cite{pambuccian2001two}: the independence model given there has cardinality~$1$.}

\begin{Proposition}\label{b8-finite-model}There exists an independence model for $\Baxiom{8}$ of cardinality~2.\end{Proposition}

\begin{proof}
Indeed, a suitable structure is already available: the countermodel $M$ for $\Baxiom{5}$ (over $\Sigma$) also falsifies $\Baxiom{8}$; since $\Baxiom{5}$ is not an axiom of $\Sigmasharpened$, $M$ works.  A counterexample to
\[
\sigmafunc{a}{b} = b \rightarrow a = b
\]
is given by $(a,b) \coloneqq (1,2)$.
In this structure, the value of $\tau$, and hence $\sigmasymbol$, is always $2$.
\end{proof}

\section{Complete independence}

%% \begin{Proposition}$\Sigmasharpened$ is independent.\end{Proposition}
%% \begin{proof}The independence proofs from~\cite{pambuccian2001two} are still valid for $\Sigmasharpened$ (putting aside the case of $\Baxiom{5}$, which, according to Theorem~\ref{b5-dependent}, is not actually independent).  If one prefers finite independence models, the theorems of Section~\ref{sec:finite-independence-models} can be used.\end{proof}

The notion of completely independent set was proposed by E.~H.~Moore~\cite{moore1910introduction}.  It is a considerably stronger property of an axiom system than the familiar notion of independence.
\begin{Definition}%[Completely independent axiom system]
An axiom system $X$ is said to be completely independent if, for all subsets $A$ of $X$, the set $A \cup \{ \phi \in X \setminus A : \neg\phi \}$ is satisfiable.
\end{Definition}
If an axiom system $X$ is completely independent then it is also dependent: for every sentence $\phi$ of $X$, we have that $X \setminus \{ \phi \} \cup \{ \neg \phi \}$ is satisfiable, or (by the completeness theorem), that $\notproves{X \setminus \{ \phi \}}{\phi}$.  When an axiom system is completely independent, no Boolean combination of its axioms can be proved from the other axioms.
\begin{Theorem}
$\Sigma$ is completely independent.
\end{Theorem}
\begin{proof}
\par{Since $\Sigma$ has $8$ axioms, by following the definition of complete independence one sees that there are $2^{8}$ sets of formulas to check for satisfiability.  (Such an enumeration of cases is best executed mechanically rather than by hand; we were assisted by the {\tipi} program~\cite{alama2012tipi}.) But for all cases, very small finite models can be found with the help of a finite model-finder for first-order classical logic (e.g., \paradox~\cite{claessen2003new}).  For lack of space we do not present all the models here.}
\end{proof}
\begin{Theorem}$\Sigmasharpened$ is completely independent.\end{Theorem}
\begin{proof}
\par{As with $\Sigma$, $\Sigmasharpened$ has $8$ axioms, so again one has $2^{8}$ sets of formulas to check for satisfiability.  Except for two cases, very small finite models can be found almost immediately.  The only cases---both involving $\Baxiom{4}$---that cannot be immediately dispensed with are:}
\begin{enumerate}
\item {$\Sigmasharpened \setminus \{ \Baxiom{4} \} \cup \{ \neg\Baxiom{4} \}$, and}
\item {$\Sigmasharpened \setminus \{ \Baxiom{4} , \Baxiom{7} \} \cup \{ \neg\Baxiom{4} , \neg\Baxiom{7} \}$.}
\end{enumerate}
\par{Case (1): The satisfiability of $\Sigmasharpened \setminus \{ \Baxiom{4} \} \cup \{ \neg\Baxiom{4} \}$ is, by the completeness theorem, the same thing as the independence of $\Baxiom{4}$.  The proof in~\cite{pambuccian2001two} works.  Recall that the smallest independence model for this axiom (27) is much larger than the other independence models (which are all size 3 or less).}
\par{Case (2): $\Sigmasharpened \setminus \{ \Baxiom{4} , \Baxiom{7} \} \cup \{ \neg\Baxiom{4} , \neg\Baxiom{7} \}$ is satisfiable.  Take a model $M$ of $\subtract{\Sigmasharpened}{\Baxiom{4}}$ in which $\Baxiom{4}$ is false.  Since $\Baxiom{4}$ fails, there exists points $a$, $b$, $c$, and $d$ in $M$ such that $\collinear{a}{b}{c}$ but $\neg\collinear{d}{\taufunc{a}{b}{d}}{\taufunc{a}{c}{d}}$.  An appropriate model is obtained from $M$ by changing $M$'s interpretation of $a_{0}$, $a_{1}$, and $a_{2}$ to, respectively, $a$, $\taufunc{a}{b}{d}$, and $\taufunc{a}{c}{d}$.}
\end{proof}
\par{The treatment of case (2) in the preceding proof might be regarded as somewhat odd.  The model at work there did contain triangles (that is, it was ``nonlinear''), yet it falsified $\Baxiom{7}$.  But if $\Baxiom{7}$ is false, shouldn't the model be ``linear''?}
\par{It is worth noting that $\Sigmasharpened$ is just barely completely independent; a seemingly innocent change to one it is axioms destroys complete independence.  Consider the existential generalization of $\Baxiom{7}$}
\begin{description}
\item[B7$^{\prime}$] {$\exists a,b,c \left [ \neg\collinear{a}{b}{c} \right ]$}
\end{description}
\par{and let $\Sigmasharpened^{\prime}$ be $\add{\subtract{\Sigmasharpened}{\Baxiom{7}}}{\Baxiom{7^{\prime}}}$.}
\par{Intuitively, $\Baxiom{7^{\prime}}$ says the same thing as $\Baxiom{7}$.  Every model of $\Sigmasharpened$ is, by ignoring the interpretations of $a_{0}$, $a_{1}$, and $a_{2}$, a model of $\Sigmasharpened^{\prime}$.  And every model of $\Sigmasharpened^{\prime}$ can be extended to a model of $\Sigmasharpened$ by choosing, for the interpretation of $a_{0}$, $a_{1}$, and $a_{2}$, any witness to the truth of $\Baxiom{7^{\prime}}$.  Nonetheless, the two theories are subtly different:}
\begin{Proposition}\label{not-completely-independent}$\Sigmasharpened^{\prime}$ is independent but not completely independent.\end{Proposition}
\begin{proof}
\par{The independence proofs for $\Sigmasharpened$ are easily adapted to $\Sigmasharpened^{\prime}$.  As for complete independence, note that $\Baxiom{7^{\prime}}$ is true in every model of $\add{\subtract{\Sigmasharpened^{\prime}}{\Baxiom{4}}}{\neg\Baxiom{4}}$.  Thus, $\Sigmasharpened^{\prime} \setminus \{ \Baxiom{4}, \Baxiom{7^{\prime}} \} \cup \{ \neg\Baxiom{4}, \neg\Baxiom{7^{\prime}} \}$ is unsatisfiable.}
\end{proof}
\par{In other words, if $\Baxiom{4}$ is false, then there is a triangle, i.e., $\Baxiom{7^{\prime}}$ holds.  This illustrates a relationship among the axioms of $\Sigmasharpened^{\prime}$ that the notion of complete independence is intended to rule out.}
\par{$\Sigmasharpened^{\prime}$ is quite far from being completely independent.  Although $\Baxiom{7^{\prime}}$ cannot be proved from the other axioms, with seven exceptions (see Table~\ref{tab:compatible-with-neg-b7-prime}), if any of $\Sigmasharpened^{\prime}$'s axiom is negated, $\Baxiom{7^{\prime}}$ becomes provable; that is, the other $2^{8} - 7$ Boolean combinations are incompatible with $\neg\Baxiom{7^{\prime}}$.  In a rough sense, then, $\Baxiom{7^{\prime}}$ is ``almost'' a theorem of $\subtract{\Sigmasharpened^{\prime}}{\Baxiom{7^{\prime}}}$.}
\begin{center}
\begin{table}[h]
\begin{tabular}{llllllll}
\hline
$\Aaxiom{3}$ & $\Baxiom{1}$ & $\Baxiom{2}$ & $\Baxiom{3}$ & $\Baxiom{4}$ & $\Baxiom{5}$ & $\Baxiom{6}$ & $\Baxiom{8}$\\
\hline
$+$ & $+$ & $+$ & $-$ & $+$ & $+$ & $-$ & $-$\\
$+$ & $+$ & $+$ & $+$ & $+$ & $+$ & $-$ & $-$\\
$+$ & $+$ & $+$ & $-$ & $+$ & $+$ & $+$ & $-$\\
$+$ & $+$ & $+$ & $-$ & $+$ & $+$ & $-$ & $+$\\
$+$ & $+$ & $+$ & $+$ & $+$ & $+$ & $+$ & $-$\\
$+$ & $+$ & $+$ & $+$ & $+$ & $+$ & $-$ & $+$\\
$+$ & $+$ & $+$ & $-$ & $+$ & $+$ & $+$ & $+$\\
\end{tabular}\caption{\label{tab:compatible-with-neg-b7-prime}Boolean combinations of axioms of $\Sigmasharpened$ compatible with $\neg\Baxiom{7^{\prime}}$.}
\end{table}
\end{center}
\par{The difference between $\Baxiom{7}$ and $\Baxiom{7^{\prime}}$ is now clear.  $\Baxiom{7}$ makes an assertion about three specific (though undetermined) points which are not mentioned anywhere else in the axioms and are thus ``semantically inert''.  By contrast, $\Baxiom{7^{\prime}}$ is a purely existential sentence that can ``interact'' with the other axioms (specifically, $\Baxiom{4}$).}
\par{Similarly, in the foundations of logic, a result similar to Proposition~\ref{not-completely-independent} was discovered by Dines~\cite{dines1915complete}: among several axioms, only one (also having the flavor of a minimal-cardinality principle) was an obstacle to complete independence.}
\par{Interestingly, the axiom system $\Gamma$ from which $\Sigma$ is derived has many dependent axioms.  $\Gamma$ is, moreover, very far from being completely independent.  Thus $\Sigma$ is, from a certain methodological standpoint, to be preferred to $\Gamma$.}
\bibliographystyle{amsplain}
\bibliography{../hyperbolic}
\end{article}
\end{document}